\newtheorem{lemma}{Lemma}[section]
\newtheorem{prop}[lemma]{Proposition}
\newtheorem{theorem}[lemma]{Theorem}
\newtheorem{cor}[lemma]{Corollary}
\newtheorem{rem}[lemma]{Remark}

\newcommand{\re}{\begin{rem}\rm}
\newcommand{\mar}{\end{rem}}

\newcommand{\kla}{\left ( }
\newcommand{\mer}{\right ) }

\newcommand{\for}{\begin{eqnarray*}}
\newcommand{\mel}{\end{eqnarray*}}

\newcommand{\mitt}{\left | { \atop } \right.}
\newcommand{\kl}{\pl \le \pl}
\newcommand{\gl}{\pl \ge \pl}
\newcommand{\kll}{\p \le \p}

\newcommand{\nach}{\rightarrow}

\newcommand{\nz}{{\rm  I\! N}}
\newcommand{\nen}{n \in \nz}
\newcommand{\ken}{k \in \nz}
\newcommand{\rz}{{\rm  I\! R}}

\newcommand{\p}{\hspace{.05cm}}
\newcommand{\pl}{\hspace{.1cm}}
\newcommand{\pll}{\hspace{.3cm}}
\newcommand{\pla}{\hspace{1.5cm}}
\newcommand{\hz}{\vspace{0.5cm}}

\newcommand{\qed}{\hspace*{\fill}$\Box$\hz\pagebreak[1]}

\newcommand{\Om}{\Omega}

\renewcommand{\a}{\alpha}
\newcommand{\si}{\sigma}

\newcommand{\eps}{\varepsilon}

\newcommand{\ce}{{\tt c}_o}

\renewcommand{\L}{{\cal L}}
\newcommand{\F}{{\cal F}}
\newcommand{\A}{{\cal A}}
\renewcommand{\P}{\Pi}
\newcommand{\PP}{{\rm I\! P}}

\newcommand{\G}{\Gamma}
\newcommand{\g}{\gamma}
\newcommand{\I}{{\cal I}}

\newcommand{\noo}{\left \|}
\newcommand{\rrm}{\right \|}

\newcommand{\bet}{\left |}
\newcommand{\rag}{\right |}

\newcommand{\intt}{\int\limits}
\newcommand{\summ}{\sum\limits}

\renewcommand{\i}{\subset}
\newcommand{\dd}{\rm I\!D}
\newcommand{\pp}{\rm I\!P}

\documentstyle[11pt]{article}
\begin{document}

\title{ \bf Hyperplane conjecture for quotient spaces of $L_p$}

\author{ Marius Junge }

\date{}

\maketitle
\begin{abstract}
\hspace{-1.9em} We give a positive solution for the hyperplane
conjecture
of quotient spaces F of $L_p$, where $1<p\kll\infty$.
\[ vol(B_F)^{\frac{n-1}{n}} \kl c_0 \pl p' \pl \sup_{H \p hyperplane}
vol(B_F\cap H)  \pl.\]
This result is extended to Banach lattices which does not contain
$\ell_1^n$'s
uniformly. Our main tools are tensor products and minimal volume ratio
with respect to $L_p$-sections.
\end{abstract}

\section*{Introduction:}
\newtheorem{heorem}{Theorem}
An open problem in the theory of convex sets is the so called \hz

{\bf Hyperplane problem:} {\it Does there exist a universal constant $c
> 0$ such that for
all $\nen$ and all convex, symmetric bodies $K \i \rz^n$ one has}
\[ |K|^{\frac{n-1}{n}} \kl c \pl \sup_{H \p hyperplane} |K\cap H|  \pl
?\]

For some classes of convex bodies this problem has a positive
solution.
For example, for convex bodies with unconditional basis
a positive solution was first given by Bourgain.  He also proved
\cite{BOU}
\[ |K|^{\frac{n-1}{n}} \kl c_0 \pl n^{\frac{1}{4}} \pl (1+\ln n) \pl
\sup_{H \p hyperplane} |K\cap H|  \pl ,\]
which is still the best known estimate for arbitrary convex bodies.
 Another class consists of convex bodies with small
volume ratio with respect to the ellipsoid of minimal volume. This
includes
the class of zonoids. K. Ball \cite{BA} solved the problem for the
duals of zonoids, i.e. unit balls of subspaces
of an $L_1$-space, briefly $L_1$-sections.

\begin{heorem}[K. Ball] For a convex, symmetric body $K \subset  \rz^n$
one has
\[ |K|^{\frac{n-1}{n}} \kl 2 \pl \pl \inf \left \{ \kla \frac{|S|}{|K|}
\mer^{\frac{1}{n}}
\mitt K\subset S,\pl S \pl L_1\mbox{-section} \right \} \pl
 \sup_{H \p hyperplane} |K\cap H|  \pl .\]
\end{heorem}

This theorem implies all positive solutions listed before (except
Bourgain's $n^{\frac{1}{4}}$ estimate).
A further application of Ball's theorem for the hyperplane problem
can be deduced from the following
\pagebreak

\begin{samepage}
\begin{heorem} Let $Y$ be a Banach space with the Gordon-Lewis
property. Then the following
conditions are equivalent.
\begin{enumerate}
\item[i)] Y does not contain $\ell_{\infty}^n$'s uniformly.
\item[ii)] There exists a constant $c>0$ such that for all $\nen$ and
all $n$-dimensional
subspace $F\subset X$ there is a $L_1$-section $S\subset F$ with
\[ \kla \frac{|S|}{|B_F|} \mer^{\frac{1}{n}} \kl c \pl .\]
\end{enumerate}
\end{heorem}
\end{samepage}

The Gordon-Lewis property was introduced in connection with a problem
of Grothendieck.
A Banach space $X$ has the Gordon-Lewis property if every absolutely
$1$-summing operator acting on $X$
factors through $L_1$. This is an operator ideal property which
typically holds
in Banach lattices. Unfortunately, Gordon and
Lewis \cite{GL} discovered spaces without this property. A combination
of theorem 1 and theorem 2 gives a positive solution of the hyperplane
conjecture
for subspaces of a Banach lattice with finite cotype. This was proved
independently by
J. Zinn (still unpublished). Since we are using the Gordon-Lewis
property it is not
surprising that the proof of theorem 2 is based on the theory of
absolutely summing operators.
But in this framework one can also replace 1-summing operators
by p-summing operators. Thereby one obtains
the theory of minimal $L_p$-sections, i.e. affine images of finite
dimensional sections of the unit ball of $L_p$,
containing a certain convex body. In this way a connection to the
case $p=2$ which is about minimal ellipsoids containing a convex body
is established. This
was discussed intensively in the literature under the name of weak
type 2 spaces. We will prove

\begin{heorem} Let $Y$ be a Banach space with the Gordon-Lewis property
and
not containing $\ell_{\infty}^n$'s uniformly. For a subspace $X
\subset  Y$ the
following conditions are equivalent.
\begin{enumerate}
\item[i)] X does not contain $\ell_1^n$'s uniformly.
\item[ii)] There exists $1<p\kll s<\infty$ and a constant $c>0$ such
that
for all $T: L_p \nach L_s$ one has
\[ \noo T\otimes Id_X : L_p(X) \nach L_s(X) \rrm \kl c \pl \noo T \rrm
\pl .\]
\item[iii)] There exists a $1 < p \kll 2$ and a constant $c_p>0$ such
that for all $\nen$ and all $n$-dimensional
subspace $F\subset X$ there is a $L_p$-section $S_p\subset F$ with
\[ \kla \frac{|S_p|}{|B_F|} \mer^{\frac{1}{n}} \kl c \pl .\]
\end{enumerate}
\end{heorem}

Analyzing the proof of theorem 3 it turns out that the type index of
$X$ coincides
with the supremum over all $p$ such that $X$ has the $L_p$-section
property. Examples show
that a type $p$ condition in Banach lattices does not imply the
$L_p$-section
property.
On the other hand $p$-convex Banach lattices with finite
cotype have the $L_p$-section property. Condition $ii)$ is related to
Bourgain's Hausdorff-Young inequalities for spaces which does not
contain $\ell_1^n$'s
uniformly. In the presence of the Gordon-Lewis property the proof is
considerably
easier and extends to arbitrary operators. But this phenomena also
indicates
a limitation of the method which is motivated by results of Pisier
\cite{PS3,PS4}.
Namely, the Schatten classes ${\cal S}_p$, $p\neq 2$, can not be
embedded in a quotient
of a subspace $X \subset Y$ which satisfies one of the above
conditions.
But on the other hand condition $ii)$ is extremely useful for the
solution of the
hyperplane conjecture.
The main idea consists in comparing  gaussian variables and coordinate
functionals
on convex bodies.

\begin{heorem} Let $Y$ be a Banach space with the Gordon-Lewis property
and not containing
$\ell_1^n$'s uniformly. Then there exists a constant $c_Y>0$ such that
for all
$n$-dimensional quotients of subspaces $F$ of $Y$ one has
\[ |B_F|^{\frac{n-1}{n}} \kl c_Y \pl \sup_{H \p hyperplane} |B_F\cap
H|  \pl.\]
In particular, the hyperplane conjecture is uniformly satisfied for
quotients
of $L_p$ for $1 < p \kl \infty$.
\end{heorem}

In section 1 we develop the tensor product techniques used for the
geometric
applications in section 2. Parts of the results are contained in the
author's
PHD-Thesis. An investigation of $L_q$-sections contained in convex
bodies
is planed in a further publication.
\section*{Preliminaries}
We are only dealing with Banach spaces over the scalar field $\rz$
of real numbers. In the text standard Banach space notation will be
used. In particular, $c_0$,
$c_1$,.. always denote universal constants. Banach spaces will be
denoted
by $E,F,Y,X,X_0,X_1,..$. The symbols $E,F$ are reserved for finite
dimensional spaces.
Given a closed subspace X of Y there is a natural injection
$\iota_X : X \nach Y$, $x \mapsto x$. We use the same notation
$\iota_X$ for
the isometric embedding of X in it's bidual. The unit ball of a Banach
space X is
denoted by $B_X$. In contrast to this we denote by $B_p^n$  ($1\le p\le
\infty$, $\nen$) the unit ball of the
space $\ell_p^n$. This space as well as  $\ell_p$, $L_p$, $\ce$ and the
vector
valued spaces $L_p(\Om,\mu;X)$, or briefly $L_p(X)$, are defined
in the usual way, where $(\Om,\mu)$ is a measure space. In the
following $p'$
denotes the conjugate index p, i.e. $\frac{1}{p} +\frac{1}{p'}={1}$.

A standard reference on operator ideals is the monograph of Pietsch
\cite{PIE}. The ideal
of all linear bounded, finite rank operators is denoted by $\L$, $\F$,
respectively.
For a Banach ideal $(\A,\alpha)$ the component $\A^*(X,Y)$ of the
conjugate ideal $(\A^*,\a^*)$ is the
class of all operators $T \in \L (X,Y)$ such that
\[     \a^*(T) \pl :=\pl \sup \left \{ |tr TS | \mitt S\in \F(Y,X) , \a
(S) \kll 1,\pl  \right \}  \pl < \pl \infty\pl .\]
The component $\A^d(X,Y)$ consists of all operators $T\in \L(X,Y)$ such
that $T^*\in \A (Y^*,X^*)$.
Equipped with the norm $\a^d(T)\p :=\p \a(T^*)$ the pair $(\A^d,\a^d)$
is again a Banach ideal.
For $1 \kll p \kll \infty$ an operator $T\in \L(X,Y)$ is said to be
(absolutely) p-summing
$(\p T \in \P_p(X,Y)\p)$ if there is a constant $c>0$ such that for all
$\nen$, $(x_k)_1^n
\subset X$
\[ \kla \summ_1^n \noo Tx_k \rrm^p \mer^{\frac{1}{p}} \kl
   c \pl \sup_{x\in B_{X^*}} \kla \summ_1^n |\langle x_k, x^* \rangle
   |^p \mer^{\frac{1}{p}}
\pl .\]
We denote by $\pi_p(T)\p:=\p \inf\{c\}$, where the infimum is taken
over all c satisfying
the above
inequality. An operator $T \in \L (X,Y)$ is $p$-integral $(\p\! T \in
\I_p (X,Y)\p\!)$, if there is a
factorization $\iota_X T\p\!=\p\! S I R$, where $R\in \L
(X,L_{\infty}(\Om,\mu))$,  $(\Om,\mu)$ a probability space,
$I \in \L (L_{\infty}(\Om,\mu), L_p(\Om,\mu))$ the formal identity and
$S \in \L (L_p(\Om,\mu), Y^{**})$.
The $p$-integral norm $\iota_p(T)$ is defined as $\inf\{ \noo S\rrm
\noo R \rrm\}$, where the infimum
is taken over all such factorizations. Let us note that if one of the
spaces $X,Y$ is finite dimensional
one has, see \cite{PIE}
\[ \P_p^*(X,Y) \pl = \pl  \I_{p'} (X,Y)\pll \mbox{and} \pll
 \I_p^* (X,Y) \pl =\pl \P_{p'} (X,Y)  \pl.\]
An operator is called p-factorable $(\p T\in \G_p (X,Y)\p)$ if there
are a measure space $(\Om,\mu)$
and operators $R\in \L (X,L_p(\Om,\mu))$, $S\in
\L(L_p(\Om,\mu),Y^{**})$ such that
$\iota_YT \p=\p SR$. Here $\iota_Y : Y\nach Y^{**}$ denotes the
canonical
embedding from Y in it's bidual. It is well known that $\G_p^d \p=\p
\G_{p'}$.

In the following $(g_k)_{\ken}$ denotes a sequence of independent,
normalized gaussian
variables on a probability space $(\Om,\PP)$. With this notion we
define for $u\in \L(\ell_2,X)$
the ideal norm
\[ \ell(u) \pl :=\pl \sup_{\nen} \noo \summ_1^n g_k u(e_k)
\rrm_{L_2(X)}\pl ,\]
where $(e_k)_{\ken}$ denotes the unit vectors in $\ell_2$ (or $\rz^n$).
In
this context Kahane's inequality is of particular interest. There
exists an absolute
constant such that for all $1 \kll p < \infty$ one has
\[ \kla \intt_{\Om} \noo \summ_1^n g_k x_k \rrm^p d\PP
\mer^{\frac{1}{p}}
\kl c_0 \pl\sqrt{p}\pl \intt_{\Om} \noo \summ_1^n g_k x_k \rrm d\PP \pl
.\]
Finally, we define a volume number $v_n(T)$ of an operator $T\in
\L(X,Y)$. This notion
is helpful to connect volume properties of Banach spaces with the
theory of operator
ideals. A closely related notion was introduced in Mascioni
\cite{MAS}.
\[ v_n(T) \pl :=\pl \sup \left \{ \kla \frac{|T(B_E)|}{|B_F|}
\mer^{\frac{1}{n}} \mitt
 E \subset X, \pl T(E)\subset F\subset Y,\pl {\rm dim}E\p=\p{\rm dim}F
 \p=\p n \right \}\pl .\]
Here and in the following $|\pl|$ denotes the translation invariant
Lebesgue measure.
If we consider the Lebesgue measure of $k$-dimensional sections of a
convex body
this will be denoted by $|\pl|_k$.
The following multiplication formula for $S\in \L(X_1,X)$ is completely
elementary
\[     v_n(TS) \kl v_n(T)\pl v_n(S)\pl  .\]
Certainly, equality holds if X is $n$ dimensional.

A Banach lattice is a Banach space with an order satisfying the same
properties as
a function spaces. For a precise definition of this and for further
notations see
\cite{LTII}. For $1 \kll p < \infty$ a Banach lattice $Y$ is said to be
p-convex, resp. p-concave if there
exists a constant $c>0$  such that for all $\nen$, $(x_k)_1^n \subset
Y$
\[ \noo \kla \summ_1^n |x_k|^p \mer^{\frac{1}{p}} \rrm \kl c \pl
	\kla \summ_1^n \noo x_k \rrm^p \mer^{\frac{1}{p}} \p,\pll
   \kla \summ_1^n \noo x_k \rrm^p \mer^{\frac{1}{p}} \kl
   c\pl \noo \kla \summ_1^n |x_k|^p \mer^{\frac{1}{p}}\rrm \p, \pl {\rm
   resp.}\]
The best possible constant will be denoted by $K^p(Y)$, $K_p(Y)$
respectively.
Let us note that $Y$ is p-convex if and only if $Y^*$ is $p'$-concave.
By a characterization of Maurey a Banach lattice is q-concave if and
only if every
positive operator $T\in \L(\ell_{\infty},Y)$ is q-summing, see again
\cite{LTII}.
Closely connected with the notion of concavity is the notion of cotype
in
arbitrary Banach spaces. A Banach space Y has cotype q ($2\kll q <
\infty$) if there is a constant $c>0$,
such that for all $\nen$, $(x_k)_1^n \subset Y$
\[ \kla \summ_1^n \noo x_k \rrm^q \mer^{\frac{1}{q}} \kl c \pl
\intt_{\Om} \noo \summ_1^n
g_k x_k \rrm d\PP \pl .\]
The best constant is denoted by $C_q(Y)$. Y has finite cotype if it has
cotype q for
some $q< \infty$. We will frequently use the following fact, which is a
combination
of Maurey's theorem and Pietsch-Grothendieck's factorization
theorem, see \cite{PIE}. Let Y be a Banach space with cotype $q$, then
for all
$q<s<\infty$ and $T \in \L (\ell_{\infty},Y)$ one has \vspace{0.2cm}

{\bf $(*)$}\hfill $\iota_s(T) \kl \pi_s(T) \kl c(s,X) \pl \noo T
\rrm\pl ,$ \hspace*{\fill}\vspace{0.2cm}

where $c(s,X)$ is a constant which only depends on $q,s$ and $C_q(X)$.
In particular,
this implies \cite[20.1.16.]{PIE} that for all $1\kll p <q'$ and $v \in
\L (X,Y)$ one has\vspace{0.2cm}

{\bf $(**)$}\hfill $\pi_1(v) \kl c(p',X) \pl \pi_p(v) \pl .$
\hspace*{\fill}\vspace{0.2cm}

\setcounter{section}{0}
\section{Gordon-Lewis property and vector-valued extensions}
In this chapter we establish the connection between the Gordon-Lewis
property with additional cotype conditions and vector valued extension
of operators between $L_p$-spaces. In a way we continue the ideas
developed
in the work of Pisier \cite{PS3}. It will be distinguish
between the Gordon-Lewis property and a restricted version $gl_2$. A
Banach
space Y is said to have the Gordon-Lewis property ( GLP), if for
every absolutely 1-summing operator $v$ the operator $\iota_Yv$
admits a factorization through some $L_1$ space. More precisely, there
exists a
constant $c>0$ such that for all $v\in \Pi_1 (Y,Z)$
\[     \gamma_1(v) \kl c \pl \pi_1(v) \pl .\]
The best possible constant is denoted by $gl(Y)$. If this inequality
only holds for
operators $v\in \Pi_1 (Y,\ell_2)$ the Banach space is said to have
$gl_2$ with
constant $gl_2(X)$.
Let us note that the GLP and $gl_2$ are self dual properties.
In the next proposition we indicate how cotype and type conditions can
be used
to improve the Gordon-Lewis-property.

\begin{prop}\label{cot} Let Y be a Banach space.
\begin{enumerate}
\item The following conditions are equivalent
\begin{enumerate}
\item[1i)] Y has the Gordon-Lewis property $ (gl_2 )$ and $Y^*$ is of
finite
cotype.
\item[1ii)]
  There exists a $1<p<\infty$ and a constant $c_p$ such that for every
  absolutely 1-summing operator $v \in \Pi_1(Y,Z)$ {\rm ($v \in
  \Pi_1(Y,\ell_2)$)}
   the operator $\iota_Yv$ admits a factorization
  through an identity $I_p: L_p(\Om,\mu ) \nach L_1(\Om,\mu)$,
  $(\Om,\mu)$
  a probability space. In other terms
\[ \iota_{p'}(v^*) \kl c_p \pl \pi_1(v) \pl .\]
\end{enumerate}
\item The following conditions are equivalent
\begin{enumerate}
\item[2i)] Y has the Gordon-Lewis property $ (gl_2 )$ and $Y$ as well
as $Y^*$ is of finite
cotype.
\item[2ii)]
  There exists $1<p,r<\infty$ and a constant $c_{pr}$ such that for
  every
  absolutely r-summing operator $v \in \Pi_r(Y,Z)$ {\rm ($v \in
  \Pi_r(Y,\ell_2)$)}
 the operator $\iota_Yv$ admits a factorization
  through an identity $I_p: L_p(\Om,\mu ) \nach L_1(\Om,\mu)$,
  $(\Om,\mu)$
  a probability space. In other terms
\[ \iota_{p'}(v^*) \kl c_p \pl \pi_r(v) \pl .\]
\end{enumerate}
\end{enumerate}
\end{prop}

{\bf Proof:} $i) \Rightarrow ii)$ By the definition of $p'$-integral
operators
it is sufficient to prove the corresponding norm inequalities. We
assume $Y^*$
of cotype $q < \infty$ and choose $p\p=\p s'$ for some $q<s<\infty$.
Let $v\in \L(Y,Z)$, ($v\in \L(Y,\ell_2)$)
be an absolutely 1-summing operator. By definition there exists a
factorization
$\iota_Yv \p=\p SR$, $R \in \L(Y,L_1)$, $S\in \L(L_1,Z^{**})$.
Using $(*)$ in the preliminaries
we see that $R^*$ is $s$ integral and get
\[ \iota_{s}(v^*) \kl \iota_s(R^*) \pl \noo S \rrm \kl c(s,Y^*) \pl
\noo R \rrm \pl \noo S \rrm \pl .\]
Taking the infimum over all factorization we obtain $c_p \kll gl(Y) \p
C(s,Y^*)$.
If $Y$ has in addition some cotype $\bar{q}$ we can apply $(**)$ in the
preliminaries
for all $1<r< \bar{q}'$.

$ii)\Rightarrow i)$ Since the identity $I_p:  L_p(\Om,\mu ) \nach
L_1(\Om,\mu)$ trivially
factors through some $L_1$ space we only have to check the
corresponding cotype conditions.
Now $1<p\kl 2$ be given by condition $ii)$ and $(y^*_i)_1^n \subset
Y^*$.
As a consequence of Kintchine's inequality one can easily see that for
\[ v \pl :=\pl \summ_1^n y^*_i \otimes e^{ }_i \pl \in \L (Y, \ell_2^n)
\]
one has
\[    \pi_1(v)  \kl  \sqrt{\frac{\pi}{2}} \pl \intt_{\Om} \noo
\summ_1^n y^*_i g^{ }_i \rrm_{Y^*}d\PP\pl.\]
Hence we deduce from the injectivity of the $p'$-summing norm
\for
 \kla \summ_1^n \noo y^*_i \rrm^{p'} \mer^{\frac{1}{p'}} &\le&
 \pi_{p'}(v^*) \kl \iota_{p'}(v^*)
\kl c_p \pl \pi_1(v) \kl c_p \pl \sqrt{\frac{\pi}{2}}  \pl \intt_{\Om}
\noo  \summ_1^n y^*_i g^{ }_i \rrm d\PP\pl.\\
\mel
This means that $Y^*$ has cotype $p'$. If condition $2ii)$ is satisfied
we use trace duality to get for all $u \in  \L(\ell_2^n,Y)$
\[ \iota_{r'}(u) \kl c_{ps} \pl \pi_p(u^*) \pl .\]
The same argument above yields that $Y$ is of cotype $r'$.\hfill
\qed\hz

\begin{rem} \label{comute} Since the GLP and $gl_2$ are self dual
proposition \ref{cot} implies
that a Banach space has GLP {\rm ($gl_2$)} and finite cotype if and
only if
there exists an $2\kll s < \infty$ and a constant $c_s$ such that every
operator $u \in \L(X,Y)$ {\rm ($u \in \L(\ell_2,Y)$)} whose dual is
absolutely 1-summing
is even $s$-integral with
\[ \iota_s(u) \kl c_s \pl \pi_1(u^*) \pl .\]
If $Y^*$ has in addition some finite cotype this improves to
\[ \iota_s(u) \kl c_s \pl \pi_p(u^*) \]
for some $p>1$.
\end{rem}

The results in this chapter are motivated by the phenomena in Banach
lattices.
In this case we can prove a sharp formula.

\begin{prop}\label{convex}
Let $1\kll p,q \kll \infty$ and Y a Banach lattice which is $p$-convex
and $q$-concave.
Then we have for all Banach space X and $u \in \L (X,Y)$
\[ \iota_q(u) \kl K^p(Y) \pl K_q(Y) \pl \pi_p(u^*) \pl .\]
\end{prop}

{\bf Proof:} Having Maurey's proof of the local unconditional structure
of Banach lattices in mind, \cite{MAU}, there is no restriction
to assume that Y has an unconditional, normalized Basis $(x_i)_{i\in
\nz}$ with
coordinate functionals $(x_i^*)^{ }_{i\in \nz}$. For $u \in \L (X,Y)$
we define
\[ S\pl:=\pl \summ_{i \in\nz} \frac{u^*(x_i^*)}{\noo u^*(x_i^*) \rrm}
\otimes e_i
\pl \in \L (X,\ell_{\infty}) \pl,\]
where $(e_i)_{i\in \nz}$ denotes the usual unit vector basis in
$\ell_p$.
We want to show that the operator $R \p:=\p \summ_{i\in\nz} e_i \otimes
\noo u^*(x_i^*) \rrm x_i$ is a positive continuous operator from
$\ell_{\infty}$
to $Y^{**}$. Clearly, $R$ is positive. For the continuity let $\a\p=\p
(\a_i)_{i\in \nz} \p
\in \p \ell_{\infty}$ and $x^* \in Y^*$ with coefficients $\beta_i
\p:=\p\langle x, x_i\rangle$.
In this situation the diagonal operator $D_{|\beta|} \in \L (\ce, Y^*)$
defines
a positive lattice homomorphism of norm at most $1$. By assumption $Y$
is
$p$-convex, hence $Y^*$ is $p'$-concave. Therefore $D_{\beta}$ is
$p'$-summing and
by Pietsch's factorization theorem there are positive diagonal
operators
$D_{\tau} \in \L (\ce,\ell_{p'})$, $D_{\si} \in \L (\ell_{p'},Y^*)$
with $D_{|\beta|}\p=\p D_{\si}\p D_{\tau}$, $\noo D_{\si} \rrm \kll 1$
and
\[\noo \tau \rrm_{p'} \kl \pi_{p'}(D_{|\beta|})\kl K_{p'}(Y^*) \pl \noo
\p|\beta |\p\rrm_{Y^*}
\kl K_p(Y) \pl \noo x^*\rrm\pl.\]
From this we obtain with H\"{o}lder's inequality
\for
|\langle R(\a), x^*\rangle| &\le& \summ_{i \in \nz} \noo u^*(x_i^*)
\rrm
	  \pl |\beta_i|  \\
&=& \summ_{i \in \nz}  \noo u^*(x_i^*) \rrm \pl \pl \si_i\pl \tau_i
\\
&\le& \kla \summ_{i \in \nz}  \noo u^*(x_i^*)\si_i
\rrm^p\mer^{\frac{1}{p}}
 \pl\noo \tau \rrm_{p'} \\
&\le& \pi_p(u^*) \pl \noo D_{\si} \rrm \pl \noo \tau \rrm_{p'}
\\[+0.3cm]
&\le& \pi_p(u^*) \pl K^p(Y) \pl \noo x^*\rrm\pl .
\mel
Let us note that $Y^{**}$ is also a $q$-concave Banach lattice.
Therefore
Maurey's characterization implies together
with Pietsch's factorization theorem applied
the operator $R$ together with $\iota_Y u\p=\p RS$
\for
 \iota_q(u) &\le& \iota_q(R) \pl \noo S \rrm \kl
K_q(Y) \pl \noo R \rrm \pl \noo  S \rrm \kl K_q(Y) \pl K^p(Y) \pl
\pi_p(u^*)\pl .\\[-1.3cm]
\mel \qed

At the end of this chapter we show how inequalities between
summing operators and integral operators can be characterized
in terms of vector-valued extensions of operators between
$L_p$-spaces.
This is connected to Kwapien's characterization of quotients of
subspaces
of $L_p$-spaces.

\begin{prop}\label{tensor}
Let $1 \kll p < \infty$, $1 \kll s \kll \infty$, $(\A,\a)$ an operator
ideal and $T\in \L (X,Y)$.
Then the following assertions are equivalent.
\begin{enumerate}
\item[i)] There exists a constant $c_1 >0$ such that for all $R \in  \A
(X_0,X_1)$,
$u \in \I_p^d (X_1,X)$
\[ \pi_s(TuR) \kl c_1 \pl \iota_p(u^*) \pl \a (R) \pl .\]
\item[ii)] There exists a constant $c_2>0$ such that for all $R \in
\A^d(\ell_p,\ell_s)$
\[ \noo R \p \otimes \p T : \ell_p(X) \nach \ell_s(Y) \rrm \kl c_2 \pl
\a (R^*)\pl. \]
\end{enumerate}
Moreover, the best constants in $i)$ and $ii)$ coincide.
\end{prop}
\pagebreak[2]

{\bf Proof:}
$i) \Rightarrow ii)$ We consider an element $u= \summ_1^n e_k \otimes
x_k \in \ell^n_p(X)$
as an operator in $\L(\ell^n_{p'},X)$ which satisfies
$\iota_p(u^*)\kll \noo u \rrm_{\ell_p^n(X)}$. In the same
way an element $w \in \ell^n_{s'}(Y^*)$ defines a $s'$-integral
operator $w \in
\I_{s'}(Y, \ell^n_{s})$ with $\iota_{s'}(w) \kll \noo w
\rrm_{\ell_{s'}^n(Y^*)}$.
With an elementary computation for traces we immediately get
\for
|\langle R\otimes T(u),w\rangle|&=&|tr(wTuR^*)|\kl \iota_{s'}(w)\pl
\pi_s(TuR^*)\kl c_1\pl \noo w\rrm \pl\noo u\rrm\pl \a (R^*)\pl .
\mel
$ii) \Rightarrow i)$ By maximality we can assume that $ii)$ also holds
for
arbitrary $L_p$, $L_s$-spaces and for $T^{**}$ instead of $T$.
Let $(x_i)_1^n \subset E$
with $\sup_{e^*\in B_{E^*}} \summ_1^n |\langle x_i, e^*\rangle|^s \kll
1$. Then
the operator $O \p:=\p \summ_1^n e_i \otimes x_i \in \L (\ell_{s'}^n,
E)$ has norm
less than one. We choose an element $w\p=\p (y_i^*)_1^n \in
B_{\ell_{s'}^n(Y^*)}$, which corresponds
to an operator from $Y$ to $\ell_{s'}^n$, such that
\[ \kla \summ_1^n \noo TuR(x_i) \rrm^s \mer^{\frac{1}{s}}
\pl =\pl \summ_1^n \langle TuR(x_i), y_i^*\rangle\pl=\pl tr(wTuRO)\pl
.\]
By the definition of $p$-integral operators there is a factorization
$u^* \p=\p SIQ$
where $Q \in \L (X^*,L_{\infty})$, $L_{\infty}$ is defined on a
probability space, $I \in \L(L_{\infty},L_p)$ the formal identity
and $S\in \L (L_p,F^*)$ such that $\noo Q\rrm \kll 1$, $\noo S\rrm\kll
(1+\eps)\p
\iota_p(u^*)$. By approximation we can even assume that the image of
$IQT^*w^*$
is contained
in the span of a finite sequence $(\chi_{A_j})_1^m$ of mutually
disjoint
characteristic functions. Then
\[ u \pl :=\pl \summ_1^m \frac{Q^*(\chi_{A_j})}{\mu (A_j)} \otimes
\chi_{A_j}\]
is an element of norm at most $1$ in $L_p(X^{**})$. We apply $ii)$ for
the operator $\bar{R}
\p:=\p O^*R^* S \in \A^d(L_p,\ell^n_{s})$ and deduce
\for
tr(wTuRO) &=& \langle \bar{R} \otimes T^{**}(u),w\rangle
\kl \noo \bar{R}\otimes T \rrm \pl \noo u\rrm_{L_p(X^{**})} \\
&\le& c_2 \pl \a(\bar{R}^*)\kl
c_2 \pl \a (R) \pl \noo O \rrm \pl \noo S\rrm
\kl c_2 \pl \a (R) \pl (1+\eps)\pl \iota_p(u^*)\pl .
\mel
Letting $\eps$ to zero yields the assertion.\qed

\re \label{translation}
For a subspace $S_p \subset L_p$ and a Banach space X we denote by
$S_p(X)$
the closure of $\{ \summ_1^m f_i \otimes x_i | f_i \in S_p, \p x_i \in
X\}$ in
$L_p(X)$. This space consists of typical $p$-summing operators
from $X^*$ to $S_p$. In a similar way the space $Q_s(X)$ is defined for
a quotient
space $Q_s$ of $L_s$. Following the same pattern as in the proof of
proposition
\ref{tensor} it can be proved that the following assertions $iii)$ and
$iv)$ as well
as $v)$ and $vi)$ are equivalent for an operator $T \in \L (X,Y)$.

\begin{enumerate}
\item[iii)] There exists a constant $c_3 >0$ such that for all $R \in
\A (X_0,X_1)$,
$u \in \Pi_p^d (X_1,X)$
\[ \pi_s(TuR) \kl c_3 \pl \pi_p(u^*) \pl \a (R) \pl .\]
\item[iv)] There exists a constant $c_4>0$ such that for all subspaces
$S_p \subset \ell_p$,
$R \in \A^d(S_p,\ell_s)$
\[ \noo R \p \otimes \p T : S_p(X) \nach \ell_s(Y) \rrm \kl c_4 \pl \a
(R^*)\pl .\]

\item[v)] There exists a constant $c_5 >0$ such that for all $R \in  \A
(X_0,X_1)$,
$u \in \Pi_p^d (X_1,X)$
\[ \iota_s(TuR) \kl c_5 \pl \pi_p(u^*) \pl \a (R) \pl .\]
\item[vi)] There exists a constant $c_6>0$ such that for all subspaces
$S_p \subset \ell_p$, all
quotients $Q_s$ of $\ell_s$ and  $R \in \A^d(S_p,Q_s)$
\[ \noo R \p \otimes \p T : S_p(X) \nach Q_s(Y) \rrm \kl c_6 \pl \a
(R^*)\pl. \]
\end{enumerate}
\mar\hz

\re\label{krm} As a consequence of the preceding remark
\ref{translation} and
remark \ref{comute} we deduce the following characterization for Banach
spaces with GLP ($gl_2$) and non-trivial type,
namely for a Banach space $Y$ the following are equivalent.
\begin{enumerate}
\item[i)] Y has the GLP $(gl_2)$ and $Y$ as well as $Y^*$ are of finite
cotype.
\item[ii)] There exists $1 <p \kll s < \infty$ and a constant
$c>0$ such that for all
subspaces $S_p \subset L_p$, quotients $Q_s$ of $L_s$ and all $T \in \L
(S_p,Q_s)$ ($T \in \Gamma_2(S_p, Q_s)$)
\[ \noo T \p \otimes \p Id_Y : S_p(Y) \nach Q_s(Y) \rrm \kl c \pl \noo
T \rrm\pl
(\p\g_2(T)\p) \pl . \]
In particular, in this situation $Y$ is K-convex and does not contain
$\ell_1^n$'s uniformly, see
\cite{PS2}.
\end{enumerate}
\mar

{\bf Proof:} We will shortly indicate why $ii)$ implies the K-convexity
of $Y$.
Indeed we denote by
\[ P \p=\p \summ_{\nen} g_n \otimes g_n \pl,\]
the orthogonal projection onto the span of a sequence of independent
normalized gaussian variables. By Kintchine's or Kahane's inequality
$P:L_p \nach L_s$ can be factorized in the form $P \p=\p u_s
(u_{p'})^*$ where
\[ u_s \p:=\p  \summ_k e_k \otimes g_k \in \L(\ell_2,L_s) \pl. \]
Therefore, we get
\[\gamma_2(P: L_p \nach L_s) \kl c_0^2 \pl \sqrt{sp'} \pl .\]
Condition $ii)$ implies together with Kahane's inequality
\[ \noo P \otimes Id_Y: L_2(Y) \nach L_2(Y) \rrm \kl c_0^2 \pl c \pl
\gamma_2(P)
\kl c \pl c_0^4 \pl \sqrt{sp'} \pl.\]
Therefore, $Y$ is K-convex and does not contain $\ell_1^n$'s uniformly,
see \cite{PS2}. \qed
\section{Geometric applications}
\setcounter{lemma}{0}
In the following a convex body will be a convex, symmetric, compact set
$K \subset \rz^n$
with non-empty interior. By $X_K \p:=\p (\rz^n,\noo \p \rrm_K)$ we
denote the
$n$-dimensional Banach space whose unit ball is $K$. The following
lemma is well-known
and can be deduced from Pajor-Tomczak's inequality and Kahane's
inequality,
for more information and constants see \cite{PS2} and \cite{SCH}. We
want to formulate this lemma
because of the frequent use.

\begin{lemma} \label{basic}
Let $1\kll s < \infty$ then for all Banach space $X$ and $u \in
\Pi_s(\ell_2,X)$
one has
\[ \sup_{\ken}  \sqrt{k}\pl v_k(u) \kl \ell (u) \kl \sqrt{s} \pl
\pi_s(u)\pl .\]
\end{lemma}


Now we are able to prove the connection between minimal $L_p$-sections
and volume
estimates for $p$-summing operators. This is a generalization of Ball's
characterization of
the weak-right-hand Gordon-Lewis property.

\begin{prop} \label{char}
There is a constant $c_0 >0$ such that for $1\kl p< \infty$
and for all convex bodies
$K \subset\rz^n$ one has
\for
\sqrt{\frac{n}{p}} \pl \sup_{\pi_p(u^*) \le 1}  v_n(u)
&\le& \inf\left \{ \kla \frac{|S_p|}{|K|}\mer^{\frac{1}{n}} \mitt K
\subset S_p,\pl S_p \pl L_p\mbox-{section} \right \}
\kl c_0 \pl \sqrt{n} \pl \sup_{\pi_p(u^*) \le 1}  v_n(u) \pl ,
\mel
where the supremum is taken over all operators $u \in \L
(\ell_2^n,X_K)$.
\end{prop}

{\bf Proof:} Let $K \subset S$ where S is an $L_p$-section, i.e.
$S\p=\p T^{-1}(B_{L_p})$
where $T \in \L (X_K,L_p)$ is a rank n operator of norm at most $1$.
For an operator
$u\in \L (\ell_2^n, X_K)$ we consider the composition $U\p:=\p Tu$
which satisfies
$\pi_p(U^*) \kll \pi_p(u^*)$. By lemma \ref{basic} and proposition
\ref{convex}
\for
\sqrt{n}\p v_n(u) &=& \sqrt{n} \kla \frac{|U(B_2^n)|}{|T(\rz^n) \cap
B_{L_p}|} \mer^{\frac{1}{n}} \pl
\kla \frac{|S|}{|K|} \mer^{\frac{1}{n}} \kl \sqrt{p} \pl \pi_p(U) \pl
\kla \frac{|S|}{|K|} \mer^{\frac{1}{n}}\\
&\le&  \sqrt{p} \pl \pi_p(U^*) \pl \kla \frac{|S|}{|K|}
\mer^{\frac{1}{n}}
\kl  \sqrt{p} \pl \pi_p(u^*) \pl \kla \frac{|S|}{|K|}
\mer^{\frac{1}{n}}
\mel
Taking the infimum over all $L_p$-sections yields the first estimate.
For the second one
we apply Lewis lemma \cite{LEW} to find an isomorphism $u \in \L
(\ell_2^n, X_K)$ with
$\pi_p(u^*) \p=\p \iota_{p'}((u^*)^{-1})\p\!=\p\! \sqrt{n}$. By
definition there is a
factorization $u^{-1} \p\!=\p\! VIR$ where  $R \in \L (X_K,L_p)$, $I
\in \L (L_p,L_1)$
the formal identity and $V\in \L (\L_1,\ell_2^n)$ with
$\noo R \rrm \!\kll\! 1$ and $\noo V \rrm \!\kll\! \sqrt{n}$. Clearly,
$S\p\!:=\p\! R^{-1}(B_{L_p})$
is an $L_p$ section which contains $K$.
As a consequence of Grothendieck's inequality
and the fact that $\ell_2$ is of (weak) cotype $2$ we deduce
\[ \sup_{\ken} \sqrt{k} \pl v_k(V) \kl c_0 \pl \pi_2(v) \kl c_0\pl
K_G\pl \noo V \rrm \kl c_0 \pl \sqrt{n} \pl .\]
If we denote the supremum on the right hand side of our assertion by
$Sup$ we
obtain
\for
1 &=& v_n(id_{\ell_2^n}) \pl =\pl v_n(u) \pl v_n(u^{-1}) \kl v_n(u) \pl
v_n(VI) \pl v_n(R)
\kl \sqrt{n}\pl Sup \pl c_0 \pl \kla \frac{|K|}{|S|} \mer^{\frac{1}{n}}
\\[-1.3cm]
\mel\qed

\re From Kwapien's inequality between $p$-summing operators it is
evident
that the supremum on the right hand side of the proceeding proposition
is minimal
for $p\p=\p 1$. Nevertheless, random quotients of $\ell_q^n$ ($1\kll
q\kll2$) with proportional
dimension $k \p =\p \delta n$ yield examples of spaces where the
minimal
volume ratio with respect to $L_p$-sections is worst possible. This was
discovered
by K. Ball in the case $p=1$. More precisely,
for such a random quotient $Q$ one has
\[ \inf \left \{ \kla \frac{|S|}{|B_Q|}\mer^{\frac{1}{k}} \mitt B_Q
\subset S,\pl S\pl L_1\mbox{-section}\right\}\pl \sim_{c_p}
\pl k^{\frac{1}{q}-\frac{1}{2}} \pl .\]
\mar

{\bf Proof:} We will show that for a random subspace $E\subset
\ell_{q'}^n$
of dimension $k\p=\p\delta n$, $1 \kll p <\infty$
\[ k^{\frac{1}{q}-\frac{1}{2}} \kl c_0 \pl \delta^{-\frac{1}{2}} \pl
\sup \{\p v_k(w)\p  | \p  \pi_p(w) \kll \sqrt{k} \p\}\pl ,\]
where the supremum is taken over all operators $w \in \L (E,\ell_2^n)$.
Then the assertion
follows from proposition \ref{char} and the inverse of
Santa\hspace{0.01cm}\'{o}'s inequality \cite{PS2}.
By \cite{FIJ} a random subspace of $\ell_{q'}^n$ satisfies
\[ \pi_1(\iota_{q'2}^n\p\iota_E) \kl c_0 \pl \pi_2(\iota_{q'2}^n) \kl
c_0 \sqrt{n} \pl .\]
Here $\iota_{q'2}^n$ denotes the formal identity from $\ell_{q'}^n$ to
$\ell_2^n$.
On the other hand a result of Meyer and Pajor, see \cite{MEPA}, implies
\for
1  &=& \kla \frac{|E\cap B_2^n|}{|E\cap B_{q'}^n|} \frac{|E\cap
B_{q'}^n|}{|E\cap B_2^n|} \mer^{\frac{1}{k}}
\kl e \pl k^{\frac{1}{2}-\frac{1}{q}}\pl  v_k(\iota_{q'2}^n \p
\iota_E)\\
\mel
Since $\ell_q$ is a (weak) type q space and ellipsoids are
$L_p$-sections this estimate is best possible.\qed

\begin{cor}\label{ppconvex} Let $1 \kll p < \infty$ and $Y$ a
$p$-convex Banach lattice. $Y$ is of
finite cotype if and only if there is a constant $c>0$ such that for
all $\nen$,
and $n$-dimensional subspace $F$ there is a $L_p$-section $S_p \subset
F$ containing
$B_F$ with
\[ \kla \frac{|S_p|}{|B_F|} \mer^{\frac{1}{n}}   \kl c \pl .\]
\end{cor}

{\bf Proof:} If Y has finite cotype it is $s$-concave for some
$s<\infty$, see \cite{LTII}.
Using lemma \ref{basic} and proposition \ref{convex} we deduce for all
$u \in \L (\ell_2^n,F)$
\for
\sqrt{n} \pl v_n(u) &\le& \ell (\iota_F u)\kl \sqrt{s} \pl \pi_s(u)\kl
\p \sqrt{s}\pl K^p(Y) \pl K_s(Y) \pl \pi_p((\iota_F \p u)^*)\\
&\le&  \p \sqrt{s}\pl K^p(Y) \pl K_s(Y) \pl \pi_p(u^*)\pl.
\mel
By proposition \ref{char} $B_F$ is contained in a $L_p$ section with
small volume.
On the other hand a Banach lattice which is not of finite cotype
contains $\ell_{\infty}^n$'s uniformly
by Maurey/Pisier's theorem \cite{MP}. Following the proceeding remark
there are $n$-dimensional Banach spaces such that the volume ratio with
respect to
minimal $L_p$-sections is of order $n^{\frac{1}{2}}$ (as far as
$p<\infty$).
Since every $n$-dimensional Banach space can be embedded $2$-isomorphic
into some
$\ell_{\infty}^{4^n}$ we deduce that $Y$ does not have the
$L_p$-section property. \qed

The proof of corollary \ref{ppconvex} contains the main idea of this
paper.
First we establish an inequality like
\[ \pi_s(u) \kl c_{ps} \pl \pi_p(u^*) \kl c_{ps} \pl \iota_p(u^*) \]
for a Banach space $X$ with the help of the Gordon-Lewis property. By
the
first chapter this corresponds to a Fubini type inequality. (That's how
it is proved
for $L_p$.) Using Lemma \ref{basic} we use this inequality for abstract
volume estimates to deduce the $L_p$-section property and estimates
for the hyperplane problem.

\re In the case $p=1$ we can apply exactly the same proof as in
corollary \ref{ppconvex}
for a Banach space with $gl_2$, simply by replacing proposition
\ref{convex}
by remark \ref{comute}. In particular,
this yields as a proof of theorem 2. Moreover, if $Y$ is a Banach space
with $gl_2$
and $Y^*$ has cotype $q$ then $Y$ itself has finite cotype if and only
if $Y$
has the $L_p$-section property for some (for all) $1<p<q'$.
\mar

{\bf Proof of theorem 3:} Let $Y$ be a Banach space which does not
contain
$\ell_{\infty}^n$'s uniformly and has GLP. By Maurey/Pisier's theorem
\cite{MP}
$Y$ is of finite cotype. Following proposition \ref{cot} and remark
\ref{comute} there
exists an $2<s<\infty$ such that for every Banach space $X_1$ and every
operator $u \in \L(X_1,Y)$ whose dual
is absolutely 1-summing the operator u is already absolutely s-summing
with
\[ \pi_s(u) \kl \iota_s(u) \kl c_s \pl \pi_1(u^*) \pl .\]

$i) \Rightarrow ii),iii)$ Let $X \subset Y$ a space not containing
$\ell_1^n$'s uniformly. From Maurey/Pisier's
characterization of non trivial type \cite{MP} the dual space $X^*$ has
some
finite cotype $q$, say. By the injectivity of the absolutely s-summing
operators
and $(**)$ in the preliminaries
we deduce for all Banach spaces $X_1$ and all $u \in \L(X_1,X)$
\[ \pi_s(u)\kl c_s \pi_1((\iota_Y u)^*) \kl c_s\pl \pi_1(u^*) \kl
c_s\pl c(p,X) \pl \pi_p(u^*) \kl c_s\pl c(p,X) \pl \iota_p(u^*) \pl,\]
for all $1\kll p <q'$. $ii)$ follows from proposition \ref{tensor} and
$iii)$ from lemma \ref{basic} and proposition \ref{char} as in the
proof of corollary \ref{ppconvex}.

$ii) \Rightarrow i)$ has already been noticed in remark \ref{krm}.

$iii) \Rightarrow i)$ We only have to show that for $1\kll p\kll 2$
\[ \inf \left \{ \kla \frac{|S|}{|B_1^n|} \mer^{\frac{1}{n}} \mitt
B_1^n \subset S, S\pl L_p \mbox{-section} \right \} \gl c_p \pl
n^{\frac{1}{p'}} \pl ,\]
which is a well-known fact due to Maurey/Carl \cite{CA} in the theory
of (weak) type $p$ spaces.

\re
For a fixed subspace $X\subset Y$, where $Y$ has $gl_2$ and finite
cotype, the proof of theorem 3
yields the following implications given $1<r<p<q \kll2$
\[  X \mbox{ type } q \pl \Rightarrow\pl  X \pl \mbox{$L_ p$-section
property} \pl \Rightarrow \pl
 X \mbox{ weak type } p \pl \Rightarrow \pl  X \mbox{ type } r\pl. \]
This means that the type index coincides with the supremum over all p
satisfying the $L_p$-section
property. For Banach lattices X with finite cotype the situation is
slightly better ($1< p\kll 2$).
\[  X\pl p\mbox{-convex}  \pl \Rightarrow\pl  X \pl \mbox{$L_
p$-section property} \pl \Rightarrow \pl
 X \mbox{ weak type } p \pl \Rightarrow \pl  X \pl r\mbox{-convex}\pl.
 \]
Nevertheless, the Lorentz spaces $\ell_{pq}$ with $p<q<\infty$ yield
examples of Banach lattices with
type p, not having the $L_p$-section property. This can be proved using
proposition \ref{char} and\vspace{0.2cm}

{\bf $(+)$}\hfill $\pi_p^d(\iota: \ell_2^n \nach \ell_{pq}^n) \pl \sim
\pl n^{\frac{1}{p}}\pl
(1+\ln n)^{\frac{1}{q}-\frac{1}{p}}\pl .$ \hspace*{\fill}\vspace{0.2cm}
\mar

{\bf Proof of (+):} It is clearly enough to prove
\[ \pi_p(\iota: \ell_{p'1}^n \nach \ell_2^n)\pl \sim_{c_0} \pl \kla
\frac{n}{1+\ln n}\mer^{\frac{1}{p}}\pl .\]
For this we note that
\[ \noo \alpha \rrm_{p'1} \pl =\pl \sup \summ_1^n \alpha_k \p \eps_k \p
\pi(k)^{-\frac{1}{p}} \pl, \]
where the supremum is taken over all signs $\eps_k=\pm 1$ and all
permutations $\pi$ of the
number $\{1,..,n\}$. We consider the group of signs $\dd_n$, the group
of
permutations $\pp_n$ with Haar measure $\mu$, $\nu$, respectively. From
the triangle inequality in $\ell_{\frac{2}{p}}$ and Kintchine's
inequality
we deduce
\for
\kla \frac{1}{n} \summ_1^n \frac{1}{k} \mer^{\frac{1}{p}} \pl \noo
\alpha \rrm_2 &=&
 \kla \summ_1^n \kla \p \intt_{\pp_n} \pi(k)^{-1} \p  \bet
 \alpha_k\rag^p \p d\nu(\pi) \mer^{\frac{2}{p}} \mer^{\frac{1}{2}}\\
&\le& \kla\p \intt_{\pp_n} \kla \summ_1^n \bet \pi(k)^{-\frac{1}{p}}
\alpha_k \rag^2 \mer^{\frac{p}{2}} \p d\nu(\pi) \mer^{\frac{1}{p}}\\
&\le& c_0 \pl \kla \p \intt_{\pp_n} \intt_{\dd_n} \bet \summ_1^n
\pi(k)^{-\frac{1}{p}} \p \eps_k \alpha_k \rag^p \p d\mu(\eps) d\nu(\pi)
\mer^{\frac{1}{p}}\\
\mel
By the trivial part of Pietsch factorization theorem this implies
\[\pi_p(\iota:\ell_{p'1}^n \nach \ell_2^n) \kl c_0^2 \pl   \kla
\frac{n}{1+\ln n}\mer^{\frac{1}{p}}\pl .\]
This estimate is optimal since $\pi_p(\iota:\ell_{p'}^n \nach
\ell_2^n)\sim n^{\frac{1}{p}}$, see \cite{PIE}.\qed

For the hyperplane problem the isotropic position of a convex body is
of particular
interest. $K \subset \rz^n$ is said to be in isotropic position if
there exists a constant L such
that
\[ \intt_K |\langle x,\theta\rangle|^2 \p \frac{dx}{|K|} \pl =\pl
L^2\pl \noo \theta \rrm_2^2 \]
holds for all vectors $\theta\in \rz^n$. If in addition $|K|\p=\p 1$
then $L_K:= L$ is called the constant of isotropy, a detailed
discussion is contained in \cite{MIPA}. For further
applications we will give a slight generalization of Hensley's result
\cite{HEN}.

\begin{prop}\label{Hen} Let $K\subset \rz^n$ be a convex body.
\begin{enumerate}
\item[i)] If $K$ is in isotropic position and $H$ is a $n-k$
dimensional subspace of $\rz^n$ one has
\[ |K|^{\frac{n-k}{nk}} \pl \kla \frac{k+2}{k} \mer^{\frac{1}{2}}
\kl \sqrt{2\pi e}\pl L_K \pl |K\cap H|_{n-k}^{\frac{1}{k}} \pl .\]
\item[ii)] There exists an absolute constant $c>0$ and an orthogonal
matrix $O$ such that
\for
 \frac{1}{c} \p \kla  \prod_{cardA=k} |O(K)\cap H_A|_{n-k}
 \mer^{\frac{1}{\kla {n \atop k} \mer k} }
&\le& |K|^{\frac{n-k}{nk}}\\
&\le& c\p L_K \p \kla  \prod_{cardA=k} |O(K)\cap H_A|_{n-k}
\mer^{\frac{1}{\kla {n \atop k} \mer k}}\p,
\mel
where $H_A$ is the span of $\{e_i \p|\ i\notin A\}$.
\end{enumerate}
\end{prop}

{\bf Proof:}
$i)$ W.l.o.g. we can assume that $K$ has volume $1$. By \cite{J2} we
get
\for
\kla \frac{k}{k+2}\mer^{\frac{1}{2}} &\le&
\kla \intt_K \noo P_{H^{\perp}}(x)\rrm_2^2 \p dx \mer^{\frac{1}{2}} \pl
|B_2^n\cap H^{\perp}|^{\frac{1}{k}} \pl |H\cap K|^{\frac{1}{k}}\\
&=& L_K \pl \sqrt{k} \pl |B_2^k|^{\frac{1}{k}} \pl |H \cap
K|^{\frac{1}{k}}
\kl \sqrt{2\pi e} \pl L_K\pl  |H \cap K|^{\frac{1}{k}} \pl,
\mel
where $H^{\perp}$ denotes the orthogonal complement of $H$.

$ii)$ The inequality on the left side follows from a result of Meyer
\cite{MEY}.
Now let us assume that $|K|\p=\p 1$. Then there exists a selfadjoint
transformation $T : \rz^n \nach \rz^n$ with $det(T)\p=\p 1$ such that
$T(K)$
is in isotropic position. By spectral decomposition there are
orthogonal matrices
$O$, $P$ such that $T \p= \p PD_{\tau}O$ where
$D_{\tau}$ is a positive diagonal operator. Since the isotropic
position is invariant under orthogonal
transformations we can assume $T \p = \p D_{\tau}O$. By the
transformation formula we deduce
from $i)$ ($c=\sqrt{6\pi e}$)
\[ \kla c\p L_k\mer^{-k} \kl |D_{\tau} O (K) \cap H_A| \pl =\pl
\kla \prod_{i\notin A}\tau_i \mer \pl |O(K)\cap H_A| \pl. \]
If we take the product over all $A$ with cardinality $k$ we obtain
\for
\kla c\p L_k \mer^{-1}
&\le& \kla \prod_{cardA=k} \prod_{i\notin A} \mer^{\frac{1}{\kla {n
\atop k} \mer k}}
 \pl \kla \prod_{cardA=k} |O(K)\cap H_A| \mer^{\frac{1}{\kla {n \atop
 k} \mer k}}\\
&=& \kla \prod_1^n \tau_i \mer^{\frac{n-k}{nk}} \pl
 \kla \prod_{cardA=k} |O(K)\cap H_A| \mer^{\frac{1}{\kla {n \atop k}
 \mer k}}\\
&=& \kla det(T) \mer^{\frac{n-k}{nk}} \pl
 \kla \prod_{cardA=k} |O(K)\cap H_A| \mer^{\frac{1}{\kla {n \atop k}
 \mer k}}\\
\mel
But $det(T)\p=\p 1$ and we have proved the assertion. \qed

Theorem 4 in the introduction is a direct consequence of the following

\begin{theorem} Let $Y$ be a Banach space with $gl_2$ and such that $Y$
and $Y^*$ have finite cotype.
Then there exists a constant $c_Y$ such that for all quotients of a
subspace
$X$ and all convex bodies $K\subset \rz^n$ one has
\[ |K|^{\frac{n-1}{n}} \kl c_Y \pll \inf\left\{\kla
\frac{|T^{-1}(B_X)|}{|K|} \mer^{\frac{1}{n}}
\mitt T:\rz^n \nach X,\pl T(K) \subset B_X \right \} \pll
\sup_{H \p hyperplane} |K\cap H|  \pl .\]
In particular, the hyperplane conjecture is uniformly satisfied for
quotients of subspaces
of $Y$.
\end{theorem}

{\bf Proof:} Since $Y$ and $Y^{*}$ have finite cotype we can apply
remark \ref{comute}
to deduce the existence of $1 < p\kll s <\infty$ such that for all $u
\in \L(\ell_2,Y)$
\[ \iota_s(u) \kl c_Y \pl \pi_p(u^*)\pl .\]
In particular, we obtain for all Banach spaces $X_0$, all Hilbert
space factorizing operator $R \in \Gamma_2(X_0,X_1)$ and all
$u\in \L(X_1,Y)$
\[ \pi_s(uR) \kll c_Y \p \iota_p(u^*) \p \gamma_2(R)\pl . \]
From proposition \ref{tensor} we deduce  \vspace{0.2cm}

{\bf (++) } \hfill $\noo P \otimes Id_X: L_p(Y) \nach L_s(Y) \rrm \kl
c_Y\pl \gamma_2(P)$ \hspace*\fill\vspace{0.2cm}

for all $P\in \Gamma_2(L_p,L_s)$. By proposition \ref{Hen}
the assertion is proved if we give an estimate for the constant of
isotropy of a convex body $K \subset \rz^n$ with $|K|\p=\p1$.
For this we want to compare a sequence of independent, normalized
gaussian variables $(g_k)_1^n$
on $(\Om,\PP)$ with the coordinate functionals $(x_k)_1^n$ on $K$.
Therefore let
us consider
\[ P\pl := \pl \summ_1^n \frac{x_k}{L_K} \otimes g_k \in \L (L_p(K),
L_s(\Om))\pl .\]
An appropriate factorization of $P$ is given by $SR$, where
\[ R \pl := \pl \summ_1^n  \frac{x_k}{L_K} \otimes e_k \p\in \p\L
(L_p(K), \ell_2^n)\pla
\mbox{and}\pla S\pl:=\pl \summ_1^n e_k \otimes g_k\p \in \p \L
(\ell_2^n,L_s(\Om))\pl .\]
By Kahane's inequality $\noo S \rrm \kll \p\sqrt{s}$. An application of
C. Borell's lemma, see \cite[Appendix]{MIS}
and the isotropic position of $K$ yields
\for
\noo R \rrm\!\! &=&\! \!\noo R^* \rrm\!\! \pl = \pl\!\!
 \sup_{\noo \beta \rrm_2 \kll 1} \kla \intt_K \bet \langle
 \frac{\beta}{L_K},x \rangle \rag^{p'} \p dx \mer^{\frac{1}{p'}}
\pl\!\! \le \pl\!\!  c_0 \pl\! p' \!\pl \sup_{\noo \beta \rrm_2 \kll 1}
\kla \intt_K \bet \langle \frac{\beta}{L_K},x \rangle \rag^2 \p dx
\mer^{\frac{1}{2}}\!\! \kl\! \!c_0\! \pl p'\pl\! .
\mel
Hence we get $\gamma_2(P) \kll c_0 \sqrt{s}\p p'$. By (++) we have
$\noo P\otimes Id_Y \rrm \kll c_X \p c_0 \p \sqrt{s} \p p'\p$.
Clearly, such an estimate is also valid for subspaces of $Y$ and, by
duality, for quotients of subspaces.
For more precise information on the injective and surjective ideal of
operators
tensoring with $P$ see \cite{DEF}. Now let $X$ be a subspace of a
quotient of $Y$
and $T:\rz^n \nach X$ with $T(K)\subset B_X$. We define $f:= \summ_1^n
T(e_k) \otimes x_k \in L_p(K;X)$
and $S:=T^{-1}(B_X)$. With lemma \ref{basic} we conclude
\begin{samepage}
\for
2\pl L_K &\le&  \kla \frac{|S|}{|K|} \mer^{\frac{1}{n}} \pl L_K \pl
\sqrt{n} \pl
 \kla \frac{|T(B_2^n)|}{|T(\rz^n)\cap B_Y|} \mer^{\frac{1}{n}} \\
&\le&  \kla \frac{|S|}{|K|} \mer^{\frac{1}{n}} \pl \noo \summ_1^n
T(e_K)\otimes L_K\p g_k\rrm_{L_2(\Om;X)}\\
&=&  \kla \frac{|S|}{|K|} \mer^{\frac{1}{n}} \pl \noo P\otimes
Id_Y(f)\rrm_{L_s(\Om;X)}\\
&\le& c_0 \pl c_X \pl \sqrt{s} \pl p'\pl  \kla \frac{|S|}{|K|}
\mer^{\frac{1}{n}} \pl
\noo f \rrm_{L_p(K;X)} \\
&=&  c_0 \pl c_X \pl \sqrt{s} \pl p'\pl  \kla \frac{|S|}{|K|}
\mer^{\frac{1}{n}} \pl
\kla \intt_K \noo T(x) \rrm^p \p dx \mer^{\frac{1}{p}}\\
&\le&  c_0 \pl c_X \pl \sqrt{s} \pl p'\pl  \kla \frac{|S|}{|K|}
\mer^{\frac{1}{n}} \\[-1.2cm]
\mel\qed
\end{samepage}

With the second part of proposition \ref{Hen} we immediately get the
following

\begin{cor} Let $Y$ be a Banach space with $gl_2$ and such that $Y$ and
$Y^*$ have finite cotype.
Then there exists a constant $c_Y$ such that for all $n$ dimensional
quotients of subspaces $F$ one has
\[ |B_F|^{\frac{n-k}{nk}} \kl c_Y \pl \sup_{H \subset F,\p codimH=k}
|B_F\cap H|_{n-k}^{\frac{1}{k}} \pl .\]
\end{cor}


\hz
1991 Mathematics Subject Classification: 52A38, 46B45, 52A21.

Key words: Unconditional basis, volume, hyperplane conjecture.
\begin{quote}
Marius Junge

Mathematisches Seminar der Universit\H{a}t Kiel

Ludewig-Meyn-Str. 4

24098 Kiel

Germany

Email: nms06@rz.uni-kiel.d400.de
\end{quote}
\end{document}